\documentclass{amsart}
\usepackage{amsmath,amsthm}
\usepackage{amsfonts,amssymb}
\usepackage{enumerate}
\usepackage{graphicx}

 \def\al{\alpha}
 \def\be{\beta}

 \def\la{\lambda}
 \def\t{\tau}

 \newcommand{\N}{{\mathbb N}}
\newtheorem{lemma}{Lemma}

\newtheorem{theorem}[lemma]{Theorem}
\newtheorem{remark}{Remark}

\makeatletter
\@addtoreset{equation}{section}
\makeatother
\def\qed{\hfill $\Box$\smallskip}

\newcommand {\ds} {\displaystyle}

\begin{document}

\title[]
{Simple bounds for the extreme zeroes of Jacobi polynomials}

\author{Geno Nikolov}
\address{Faculty of Mathematics and Informatics, Sofia University
"St. Kliment Ohridski", 5 James Bourchier Blvd.,
1164 Sofia, Bulgaria} \email{geno@fmi.uni-sofia.bg}

\thanks{This study is financed by the European Union-NextGenerationEU,
through the National Recovery and Resilience Plan of the Republic of
Bulgaria, project No BG-RRP-2.004-0008.}

\begin{abstract}
Some new bounds for the extreme zeroes of Jacobi polynomials are
obtained with an elementary approach. A feature of these bounds is
their simple forms, which make them easy to work with. Despite their
simplicity, our lower bounds for the largest zeroes of Gegenbauer
polynomials are compatible with some of the best hitherto known
results.
\medskip

\medskip

\end{abstract}

\maketitle
\section{Introduction and statement of the results}
Zeroes of the classical orthogonal polynomials have been a topic of
permanent interest. A huge number of publications is devoted to the
study of their extreme zeroes. Many of the classical results on the
subject are collected in the Szeg\H{o} monograph \cite{Szego(1975)}.
Without any claim for completeness, we refer to \cite{ADGR(2004),
ADGR(2012), DimNik(2010), DriJor(2012), DriJor(2013), GupMul(2007),
IsmLi(1992), IsmMul(1995), Kras(2003), Kras(2003)a, Kras(2006),
GN(2019), GN(2023), NU(2018)} for some recent developments.

Our concern here is the extreme zeroes of Jacobi and Gegenbauer
polynomials. Recall that the Jacobi polynomials are orthogonal in
$[-1,1]$ with respect to the weight function
$w_{\al,\be}(x)=(1-x)^{\al}(1+x)^{\be}$, $\al,\be>-1$. The explicit
form of the $n$-th Jacobi polynomial $P_n^{(\al,\be)}$ is
$$
P_n^{(\al,\be)}(x)=\frac{1}{2^n}\,\sum_{k=0}^{n}{n+\al\choose
n-k}{n+\beta\choose k} (x-1)^{n-k}(x+1)^k
$$
(see, e.g., \cite[p. 144, eq. (2.6)]{Chihara(1978)} or \cite[p. 68,
eq. (4.3.2)]{Szego(1975)}. The Gegenbauer (also called as
ultraspherical) polynomials are orthogonal in $[-1,1]$ with respect
to the weight function $w_{\la}(x)=(1-x^2)^{\la-1/2}$, $\la>-1/2$.
The $n$-th ultraspherical polynomial $P_n^{(\la)}$ and the Jacobi
polynomial $P_n^{(\al,\al)}$ are related by the equation
$$
P_n^{(\la)}(x)={2\al\choose\al}^{-1}{n+2\al\choose\al}\,P_n^{(\al,\al)}(x),
\quad \al=\la-1/2\ne -1/2
$$
(see \cite[p. 144, eq. (2.5)]{Chihara(1978)}). Throughout the paper
we use the notation
$$
x_{1,n}(\al,\be)<x_{2,n}(\al,\be)<\cdots<x_{n,n}(\al,\be)
$$
for the zeros of the $n$-th Jacobi polynomial $P_n^{(\al,\be)}$,
$\al,\be>-1$, and the zeros of the $n$-th ultraspherical polynomial
$P_n^{(\la)}$, $\la>-1/2$ are denoted by
$$
x_{1,n}(\la)<x_{2,n}(\la)<\cdots<x_{n,n}(\la)\,.
$$
Upper bounds for the largest zero and lower bounds for the smallest
zero will be referred to as outer bounds, while by inner bounds we
mean lower bounds for the largest zero and upper bound for the
smallest zero.

For derivation of inner and outer bounds for the extreme zeros of
classical orthogonal polynomials various techniques have been
employed, among them are Sturm's comparison theorem for the zeros of
solutions of second order differential equation, the Euler-Rayleigh
method, Jensen inequalities for entire functions from the
Laguerre-P\'{o}lya class and their refinement for real-root
polynomials, formulated as a conjecture in \cite{FK(2002)} and
proved in \cite{NU(2004)}, A. Markov's and Hellmann-Feynmann's
theorems on monotone dependance of zeros of orthogonal polynomials
on a parameter, Obreshkov's extension of Descartes' rule of signs,
etc.

In \cite[Section 6.3]{Szego(1975)} sharp estimates for the zeros of
$P_n^{(\al,\be)}$ are given under the restriction $-1/2\leq
\al,\be\leq 1/2$, as well as asymptotic formulae for the zeros of
the Gegenbauer polynomials $P_n^{(\la)}$ (when $0<\la<1$), Laguerre
and Hermite polynomials through the zeros of related Bessel
functions. By using the so-called Bethe ansatz equations, Krasikov
\cite{Kras(2003), Kras(2003)a, Kras(2006)} proved two sided
estimates for the extreme zeros of Jacobi and Laguerre polynomials,
which hold uniformly with respect to all parameters involved, and
thus locate the extreme zeros of the classical orthogonal with a
high precision.

For practical purposes sometimes preference is given to bounds which
are less precise, but hold true for all values of parameters
involved and are given by simple expressions, allowing easy
manipulation with them. The aim of the present note is to prove a
simple outer bound for the extreme zeros of Jacobi polynomial
$P_n^{(\al,\be)}$ and simple inner bounds for the extreme zeros of
the Gegenbauer polynomial $P_n^{(\la)}$.

Our model for simple bounds is the outer bounds for the extreme
zeros of $P_n^{(\al,\be)}$ derived with the Newton--Raphston
iteration method:
\begin{equation}\label{e1.1}
1-x_{n,n}(\al,\be)\geq \frac{2(\al+1)}{n(n+\al+\be+1)},\qquad
1+x_{1,n}(\al,\be)\geq \frac{2(\be+1)}{n(n+\al+\be+1)}
\end{equation}
(here and in what follows, we prefer to estimate the distance
between the extreme zeros and the endpoints of $[-1,1]$). Although
extremely simple, the bounds in \eqref{e1.1} represent correctly the
behavior of the extreme zeros of $P_n^{(\al,\be)}$ when either one
of parameters $\al$ and $\be$ tends to $-1$ or to infinity, or $n$
grows.

Another example of simple outer bounds for the extreme zeros of
$P_n^{(\al,\be)}$ is due to Nevai, Erd\'{e}lyi and Magnus
\cite[Theorem 13]{NEM(1994)}. Using a variant of the Sturmian
comparison theorem, they proved that if $\al,\be\geq -\frac{1}{2}$,
then
$$
1-x_{n,n}(\al,\be)\geq\frac{2\al^2}{(2n+\al+\be+1)^2},\quad
1+x_{1,n}(\al,\be)\geq\frac{2\be^2}{(2n+\alpha+\beta+1)^2},
$$
with a further refinement for $\al,\be>0$, given by
\begin{eqnarray*}
&&1-x_{n,n}(\al,\be)\geq\frac{2\al^2}{(2n+\al)(2n+\al+2\be+2)},\\
&&1+x_{1,n}(\al,\be)\geq\frac{2\be^2}{(2n+\be)(2n+\be+2\al+2)}.
\end{eqnarray*}

Our first result is simple outer bounds for the extreme zeros of
Jacobi polynomials, which improve \eqref{e1.1} and hold for all
admissible values of $\al$ and $\be$.
\begin{theorem}\label{t1}
For all $n\geq 1$ and $\al,\be>-1$ the extreme zeros of
$P_n^{(\al,\be)}$ satisfy
\begin{equation}\label{e1.2}
1-x_{n,n}(\al,\be)\geq \frac{2(\al+1)}{n(n+\be)+\al+1},\qquad
1+x_{1,n}(\al,\be)\geq \frac{2(\be+1)}{n(n+\al)+\be+1}.
\end{equation}
\end{theorem}
Recently we applied Theorem~\ref{t1} to the study of regularity of
certain Birkhoff interpolation problems \cite{KN(2024)}.

In the opposite direction, the following simple but rather crude
inner bounds for the zeros of Jacobi polynomials are obtained with a
theorem due to Laguerre (see \cite[eq. (6.2.1)]{Szego(1975)}):
$$
1-x_{n,n}(\al,\be)\leq \frac{2(\al+1)}{2n+\al+\be},\qquad
1+x_{1,n}(\al,\be)\leq \frac{2(\be+1)}{2n+\al+\be}.
$$
Sharper simple inner bounds are obtained by Driver and Jordaan
\cite{DriJor(2012)}:
\begin{eqnarray*}
&&1-x_{n,n}(\al,\be)\leq
\frac{2(\al+1)(\al+3)}{2n(n+\al+\be+1)+(\al+1)(\al+\be+2)},\\
&&1+x_{1,n}(\al,\be)\leq
\frac{2(\be+1)(\be+3)}{2n(n+\al+\be+1)+(\be+1)(\al+\be+2)}
\end{eqnarray*}
(see also \cite[Theorems 1.4, 1.5]{GN(2019)} and \cite[Theorems 1,
2]{GN(2023)}).

Regarding estimates for the extreme zeros of Gegenbauer polynomials,
we should mention the results of Krasikov \cite{Kras(2006)}, who has
shown that
\begin{equation}\label{e1.3}
x_{n,n}(\la)=S\Big(1-\delta\frac{(1-S^2)^{2/3}}{(2R)^{1/3}S}\Big)\,,\quad
3<\delta<9,
\end{equation}
where
$$
S=\sqrt{\frac{4n(n+2\la)}{4n(n+2\la)+(2\la+1)^2}},\qquad
R=2\sqrt{n(n+2\la)\big(4n(n+2\la)+(2\la+1)^2\big)}\,.
$$
The order of the error term is
$$
\frac{(1-S^2)^{2/3}}{(2R)^{1/3}S}=
O\Big(\frac{(2\la+1)^{1/3}}{n^{2/3}(n+\la)^{4/3}}\Big)\,,
$$
meaning that \eqref{e1.3} provides a second order bound, and
therefore Krasikov's bounds seem hardly improvable.

However, as already pointed out, we emphasize on the simplicity of
the bounds for the extreme zeros. This is why we quote, for the sake
of comparison, the inner bound (\cite[Theorem~3]{GN(2023)})
\begin{equation}\label{e1.4}
1-x_{n,n}(\la)<\frac{(2\la+1)(2\la+3)(2\la+7)}
{(10\la+17)\big[n(n+2\la)+\frac{1}{8}(2\la+1)^2\big]}\,,
\end{equation}
and the outer bound \cite[Lemma~3.5]{GN(2005)}
\begin{equation}\label{e1.5}
1-x_{n,n}^2(\la)>\frac{(2\la+1)(2\la+9)}{4n(n+2\la)+(2\la+1)(
2\la+5)}\,.
\end{equation}
\begin{remark}
One can derive as a consequence from \cite[eqn. (1.4)]{DimNik(2010)}
the estimate
$$
1-x_{n,n}^2(\la)>\frac{(2\la+1)\big[(2\la+9)n+4(2\la-3)\big]}
{4(n+\la-1)\big[n(n+\la-1)+4(\la+1)\big]},
$$
which is slightly sharper though less simple than \eqref{e1.5}.
Although rather sharp for a fixed $\la$, the estimate in
\eqref{e1.4} deteriorates when $n$ is fixed and $\la$ grows.
\end{remark}

Our first inner bounds for the extreme zeros of Gegenbauer
polynomials improve upon earlier results from \cite{DriJor(2012)}
and \cite{GN(2019)}.
\begin{theorem}\label{t2}
For all $n\geq 4$ and $\la>-\frac{1}{2}$, the largest zero
$x_{nn}(\la)$ of the Gegenbauer polynomial $P_n^{(\la)}$ satisfies
\begin{equation}\label{e1.6}
1-x_{n,n}^2(\la)\leq
\frac{(2\la+1)(2\la+5)}{2n(n+2\la)+(2\la+1)(2\la+2)}\,.
\end{equation}
\end{theorem}
Theorem~\ref{t2} is obtained on the way towards the proof of the
following stronger result.
\begin{theorem}\label{t3}
For all $n\geq 5$ and $\la>-\frac{1}{2}$, the largest zero
$x_{n,n}(\la)$ of the Gegenbauer polynomial $P_n^{(\la)}$ satisfies
\begin{equation}\label{e1.7}
1-x_{n,n}^2(\la)<\frac{2(2\la+1)(2\la+7)}
{c\,n(n+2\la)+4(\la+2)(2\la+7-2c\big)}\,,
\end{equation}
where
$$
c=c(\la)=3+\sqrt{5+\frac{32}{(2\la+3)(2\la+5)}}\,.
$$
\end{theorem}
Clearly, $c(\la)\in (3+\sqrt{5},6)$, and Theorem~\ref{t3} can be
stated with the constant $3+\sqrt{5}\approx 5.236$ without essential
loss of accuracy when $\la$ is large. However, for small $\la$ the
bounds with $c(\la)$ are preferable. For instance, in the case of
Chebyshev polynomials of the first and second kind ($\la=0$ and
$\la=1$), for the quantities
$$
1-x_{nn}^2(0)=\sin^2\frac{\pi}{2n}\approx \frac{\pi^2}{4n^2}\approx
\frac{2.4674}{n^2},\quad 1-x_{nn}^2(1)=\sin^2\frac{\pi}{n}\approx
\frac{\pi^2}{n^2}\approx \frac{9.8696}{n^2}
$$
Theorem~\ref{t3} yields upper bounds approximately equal to $\,\ds
\frac{2.468774}{n^2}$ and $\ds \frac{9.941217}{n^2}$, respectively,
improving the ones in \cite[p. 1802]{DimNik(2010)}.

We conclude this section with pointing out that the ratio of the
upper and the lower bound for $1-x_{n,n}^2(\la)$ given by
\eqref{e1.7} and \eqref{e1.5}, respectively, is uniformly bounded by
$2(3-\sqrt{5})\approx 1.527864$.
\section{Proofs}
\subsection{Proof of Theorem~\ref{t1}}
According to the Enestr\"{o}m--Kakeya theorem (see, e.g.,
\cite{ASV(1979)}), the modulus of each root of a polynomial
$P(z)=a_0+a_1z+\cdots+a_n z^n$ with positive coefficients is between
$\ds \min_{1\leq k\leq n}\frac{a_{k-1}}{a_k}$ and $\ds \max_{1\leq
k\leq n}\frac{a_{k-1}}{a_k}$. The Jacobi polynomial
$P_n^{(\al,\be)}$ can be represented in the following form
$$
P_n^{(\al,\be)}(x)=\Big(\frac{x+1}{2}\Big)^{n}{n+\al \choose
n}\sum_{k=0}^{n}{n\choose k}\frac{(n+\be-k+1)_k}{(\al+1)_k}
\Big(\frac{x-1}{x+1}\Big)^{k}
$$
(cf. \cite[eq. 4.3.2]{Szego(1975)}), where $(b)_0=1$ and
$(b)_k=b(b+1)\cdots(b+k-1)$ for $k\in \N$ is the Pochhammer symbol.
We apply the Enestr\"{o}m--Kakeya theorem to the polynomial $P(z)$
with
$$
a_k={n\choose k}\frac{(n+\be-k+1)_k}{(\al+1)_k},\quad
z=\frac{x-1}{x+1}.
$$
Since $\ds
\frac{a_{k-1}}{a_{k}}=\frac{k(\al+k)}{(n+1-k)(n+\be+1-k)}$ increases
with $k$, it follows from the Enestr\"{o}m--Kakeya theorem that
$$
\Big|\frac{x_{n,n}(\al,\be)-1}{1+x_{n,n}(\al,\be)}\Big|
=\frac{1-x_{n,n}(\al,\be)}{1+x_{n,n}(\al,\be)} \geq \frac{a_0}{a_1}=
\frac{\al+1}{n(n+\be)},
$$
whence the first inequality in \eqref{e1.2} follows. Invoking again
the Enestr\"{o}m--Kakeya theorem, we conclude that
$$
\Big|\frac{x_{1,n}(\al,\be)-1}{1+x_{1,n}(\al,\be)}\Big|
=\frac{1-x_{1,n}(\al,\be)}{1+x_{1,n}(\al,\be)} \leq
\frac{a_{n-1}}{a_{n}}= \frac{n(n+\al)}{\be+1},
$$
which yields the second inequality in \eqref{e1.2}. The second
inequality in \eqref{e1.2} can be deduced also from the first one
and the well-known symmetry property of Jacobi polynomials
$P_{n}^{(\al,\be)}(-x)=(-1)^n\,P_n^{(\be,\al)}(x)$ (cf. \cite[p.
144, eq. (2.8)]{Chihara(1978)}). \qed
\subsection{Proof of Theorems~\ref{t2} and \ref{t3}}
Our proof makes use of the second order differential equation (see
\cite[p. 80, eq. (4.7.5)]{Szego(1975)}
$$
(1-x^2)y^{\prime\prime}-(2\la+1)x\,y^{\prime}+n(n+2\la)y=0,\quad
y=P_n^{(\la)}(x).
$$
By differentiating this equation we obtain ordinary differential
equations satisfied by the derivatives of $y$:
\begin{equation}\label{e2.1}
(1-x^2)y^{(q+2)}(t)-(2\la+2q+1)ty^{(q+1)}(t)+(n-q)(n+2\la+q)y^{(q)}(t)=0.
\end{equation}
We combine \eqref{e2.1} with the trivial observation that for
$\t=x_{n,n}(\la)>0$ there holds
\begin{equation}\label{e2.2}
y^{(q+2)}(\t)y^{(q+1)}(\t)>0,\quad q=0,1,\ldots,n-2.
\end{equation}

By using $y(\t)=0$, we find from \eqref{e2.1} with $q=0$
$$
y^{\prime}(\t)=\frac{1-\t^2}{(2\la+1)\t}\,y^{\prime\prime}(\t).
$$
Substituting this expression for $y^{\prime}(\t)$ in \eqref{e2.1}
with $q=1$ and $t=\t$, we arrive at the equation {\small
\begin{equation*}
\begin{split}
\frac{(1-\t^2)y^{\prime\prime\prime}(\t)}{y^{\prime\prime}(\t)}
&=\frac{1}{(2\la+1)\t}\Big\{(2\la+1)(2\la+3)
-\big[n(n+2\la)+(2\la+1)(2\la+2)\big](1-\t^2)\Big\}\\
&=:\frac{A}{(2\la+1)\t}\,.
\end{split}
\end{equation*}}
Now \eqref{e2.2} with $q=1$ implies $A>0$, resulting in the inner
bound
$$
1-x_{n,n}^2(\la)<\frac{(2\la+1)(2\la+3)}{n(n+2\la)+(2\la+1)(2\la+2)},
$$
which has been obtained by another method by Driver and Jordaan
\cite{DriJor(2012)}.

Next, we express $y^{\prime\prime}(\t)$ through
$y^{\prime\prime\prime}(\t)$, and substitute the resulting
expression in \eqref{e2.1} with $q=2$ to obtain
$$
\frac{(1-\t^2)y^{(4)}(\t)}{\t
y^{\prime\prime\prime}(\t)}=\frac{(2\la+3)B}{A},
$$
where
$$
B=(2\la+1)(2\la+5)-\big[2n(n+2\la)+(2\la+1)(2\la+2)\big](1-\t^2)\,.
$$
Assuming $n\geq 4$, we conclude from \eqref{e2.2} with $q=2$ that
$B>0$, which proves Theorem~\ref{t2}.\qed

For the proof of Theorem~\ref{t3} we need to repeat the above
procedure once again, this time assuming $n\geq 5$. Lengthy but
straightforward calculations lead to the equation
\begin{equation}\label{e2.3}
\frac{(1-\t^2)y^{(5)}(\t)}{y^{(4)}(\t)}=\frac{C}{(2\la+3)\t B},
\end{equation}
where
$$
C=C(u)=K\,u^2-L\,u+M\,,\qquad u=1-\t^2,
$$
with coefficients $K,L,M$ given by
$$
K=n(n+2\la)\big[n(n+2\la)+12\la^2+40\la+35\big]
+(2\la+1)(2\la+2)(2\la+3)(2\la+4),
$$
$$
L=(2\la+3)(2\la+5)\big[3n(n+2\la)+4(\la+2)(2\la+1)\big],
$$
$$
M=(2\la+1)(2\la+3)(2\la+5)(2\la+7)\,.
$$
We show that the discriminant $\Delta=L^2-4KM$ is positive. Indeed,
by using $(2\la+3)(2\la+5)>(2\la+1)(2\la+7)$ we obtain the
inequality
\begin{equation*}
\begin{split}
\Delta>(2\la+3)^2(2\la+5)^2\Big\{&\big[3n(n+2\la)+4(\la+2)(2\la+1)\big]^2\\
&-4n(n+2\la)\big[n(n+2\la)+12\la^2+40\la+35\big]\\
&-4(2\la+1)(2\la+2)(2\la+3)(2\la+4)\Big\},
\end{split}
\end{equation*}
which simplifies to
$$
\Delta>(2\la+3)^2(2\la+5)^2\Big[5n^2(n+2\la)^2-4(10\la+23)n(n+2\la)
-16(2\la+1)(\la+2)\Big].
$$
It is easily verified that for $n\geq 5$ the expression in the
brackets is greater than
$5\big[n(n+2\la)-8(\la+2)\big]^2=5(n-4)^2(n+2\la+4)^2$, therefore
\begin{equation}\label{e2.4}
\Delta>5(2\la+3)^2(2\la+5)^2(n-4)^2(n+2\la+4)^2=:\widetilde{\Delta},\quad
n\geq 5.
\end{equation}
From \eqref{e2.3} and \eqref{e2.2} with $q=3$ we have $C>0$, hence
$u\not\in [u_1,u_2]$, where $u_{1,2}$ are the roots of quadratic
equation $C(u)=0$, i.e.,
$$
u_{1,2}=\frac{2M}{L\pm\sqrt{\Delta}}
=\frac{2(2\la+1)(2\la+7)}{3n(n+2\la)+4(\la+2)(2\la+1)
\pm\frac{\sqrt{\Delta}}{(2\la+3)(2\la+5)}}\,.
$$
The inequality $u>u_2$ is impossible. Indeed, in view of
Theorem~\ref{t2}, we have
$$
u<\frac{(2\la+1)(2\la+5)}{2n(n+2\la)+(2\la+1)(2\la+2)}
<\frac{2(2\la+1)(2\la+7)}{3n(n+2\la)+4(\la+2)(2\la+1)}<u_2.
$$
Therefore,
\begin{equation}\label{e2.5}
u=1-x_{n,n}^2(\la)<u_1
=\frac{2(2\la+1)(2\la+7)}{3n(n+2\la)+4(\la+2)(2\la+1)
+\frac{\sqrt{\Delta}}{(2\la+3)(2\la+5)}}.
\end{equation}
Replacing $\Delta$ by its lower bound $\widetilde{\Delta}$ from
\eqref{e2.4}, we arrive at the inequality
\begin{equation*}
1-x_{n,n}^2(\la)<\frac{2(2\la+1)(2\la+7)}
{c\,n(n+2\la)+4(\la+2)(2\la+7-2c\big)}\,,\quad c=3+\sqrt{5}\,.
\end{equation*}
The estimate for $1-x_{n,n}^2(\la)$ in Theorem~\ref{t3} with the
smaller constant $c(\la)$ follows from \eqref{e2.5} in the same way,
this time replacing $\Delta$ by the more precise lower bound
\begin{equation}\label{e2.6}
\Delta> \Big(1+\frac{32}{5(2\la+3)(2\la+5)}\Big)\widetilde{\Delta},
\qquad n\geq 5.
\end{equation}
For the proof of \eqref{e2.6} we find $\Delta=(2\la+3)(2\la+5)D$,
where
\begin{equation*}
\begin{split}
D=&(20\la^2+80\la+107)n^2(n+2\la)^2
-4(2\la+1)(20\la^2+68\la+65)n(n+2\la)\\&+24(2\la+1)^2(2\la+3)(2\la+4).
\end{split}
\end{equation*}
We then verify that for $n\geq 5$
\begin{equation*}
\begin{split}
D&>(20\la^2+80\la+107)\big[n(n+2\la)-8(\la+2)\big]^2\\
 &=(2\la+3)(2\la+5)\Big(5+\frac{32}{(2\la+3)(2\la+5)}\Big)(n-4)^2(n+2\la+4)^2,
\end{split}
\end{equation*}
whence \eqref{e2.6} follows. The details are left to the reader.\qed


\begin{thebibliography}{10}
\bibitem{ASV(1979)}
N. Andersen, E.\,B. Saff, R.\,S. Varga, \textit{On the Enestr\"{o}m
--Kakeya theorem and its sharpness}, Linear Alg. Appl. \textbf{28}
(1979), 5--16.

\bibitem{ADGR(2004)}
I. Area, D.\,K. Dimitrov, E. Godoy, A. Ronveaux, \textit{Zeros of
Gegenbauer and Hermite polynomials and connection coefficients},
Math. Comp. \textbf{73} (2004), 1937--1951.

\bibitem{ADGR(2012)}
I. Area, D.\,K. Dimitrov, E. Godoy, F.\,R. Rafaeli,
\textit{Inequalities for zeros of Jacobi polynomials via
Obreshkoff's theorem}, Math. Comp. \textbf{81} (2012), 991--1012.

\bibitem{Chihara(1978)}
T.\,S. Chihara, \textit{An introduction to orthogonal polynomials},
Gordon and Breach Sci. Pub., 1978.

\bibitem{DimNik(2010)}
D.\,K. Dimitrov, G.\,P. Nikolov, \textit{Sharp bounds for the
extreme zeros of classical orthogonal polynomials}, J. Approx.
Theory \textbf{162} (2010), 1793--1804.

\bibitem{DJJ(2011)}
K. Driver, A. Jooste, K. Jordaan, \textit{Stiltjes interlasing of
zeros of Jacobi polynomials from different sequences}, Ellectronic
Trans. Numer. Anal. \textbf{38} (2011), 317--326.

\bibitem{DriJor(2012)}
K. Driver, K. Jordaan, \textit{Bounds for extreme zeros of some
classical orthogonal polynomials}, J. Approx. Theory \textbf{164}
(2012), 1200--1204.

\bibitem{DriJor(2013)}
K. Driver, K. Jordaan, \textit{Inequalities for extreme zeros of
some classical orthogonal and $q$-orthogonal polynomials}, Math.
Model. Nat. Phenom.  \textbf{8}(1) (2013), 48--59.

\bibitem{FK(2002)}
W.\,H. Foster and I. Krasikov, \textit{Inequalities for real-root
polynomials and entire functions}, Adv. Appl. Math. {\bf 29} (2002),
102--114.

\bibitem{GupMul(2007)}
D.\,P. Gupta, M.\,E. Muldoon, \textit{Inequalities for the smallest
zeros of Laguerre polynomials and their $q$-analogues}, J. Ineq.
Pure Appl. Math. \textbf{8}(1) (2007), Article 24.

\bibitem{IsmLi(1992)}
M.\,E.\,H. Ismail, X. Li, \textit{Bounds on the extreme zeros of
orthogonal polynomials}, Proc. Amer. Math. Soc. \textbf{115} (1992),
131--140.

\bibitem{IsmMul(1995)}
M.\,E.\,H. Ismail, M.\,E. Muldoon, \textit{Bounds for the small real
and purelyimaginary zeros of Bessel and related functions}, Met.
Appl. Math. Appl. \textbf{2} (1995), 1--21.

\bibitem{KN(2024)}
B. Konstantinova, G. Nikolov, \textit{Two sequences of regular
three-row almost Hermitian incidence matrices}, to appear in:
Constructive Theory of Functions (B. Draganov et al.,Eds.), 2024.

\bibitem{Kras(2003)a}
I. Krasikov, \textit{On the zeros of polynomials and allied
functions satisfying second order differential equations}, East J.
Approx. \textbf{9} (2003), 51--65.

\bibitem{Kras(2003)}
I. Krasikov, \textit{Bounds for zeros of the Laguerre polynomials},
J. Approx. Theory \textbf{121} (2003), 287--291.

\bibitem{Kras(2006)}
I. Krasikov, \textit{On extreme zeros of classical orthogonal
polynomials}, J. Comp. Appl. Math.\textbf{193} (2006), 168--182.

\bibitem{NEM(1994)}
P. Nevai, T. Erd\'{e}lyi, A. Magnus, \textit{Generalized Jacobi
weights, Christoffel functions, and Jacobi polynomials}, SIAM J.
Math. Anal. \textbf{25}(2) (1994), 602--614.

\bibitem{GN(2005)}
G. Nikolov, \textit{Inequalities of Duffin--Schaeffer type, II},
East J. Approx. \textbf{11} (2005), 147--168.

\bibitem{GN(2019)}
G. Nikolov, \textit{New bounds for the extreme zeros of Jacobi
polynomials}, Proc. Amer. Math. Soc. \textbf{147}(4) (2019),
1541--1550.

\bibitem{GN(2023)}
G. Nikolov, \textit{On the extreme zeros of Jacobi polynomials}, in:
Numerical Methods and Applications 2022 (I. Georgiev et. al, Eds.)),
LNCS 13858, Springer Nature Switzerland 2023, pp. 246--257.

\bibitem{NU(2004)}
G. Nikolov, R. Uluchev, \textit{Inequalities for real root
polynomials. Proof of a conjecture of Foster and Krasikov}, in:
``Approximation Theory: A Volume Dedicated to Borislav Bojanov''
(D.\,K. Dimitrov, G. Nikolov, and R. Uluchev, Eds.), Marin Drinov
Academic Publishing House, Sofia 2004, pp. 201--216.

\bibitem{NU(2018)}
G. Nikolov, R. Uluchev, \textit{Bounds for the extreme zeros of
Laguerre polynomials}, in: Numerical Methods and Applications 2018
(G. Nikolov et. al, Eds.)), LNCS 11189, Springer Nature Switzerland
2019, pp. 243--250.

\bibitem{Szego(1975)}
G. Szeg\H{o}, \textit{Orthogonal Polynomials}, 4th edn. AMS
Colloquium Publications, Providence, RI (1975).

\end{thebibliography}
\end{document}